\documentclass[a4paper,reqno,10pt]{amsart}
\usepackage{amssymb}
\usepackage{amsmath}
\usepackage{amsthm}
\usepackage{graphicx}
\usepackage{color}
\usepackage{epsfig}
\usepackage{amsfonts}
\usepackage{hyperref}
\usepackage{mathtools}
\usepackage{bm}

\newcommand{\R}{\mathbb R}
\newcommand{\N}{\mathbb N}

 \allowdisplaybreaks \numberwithin{equation}{section}
\begingroup
\theoremstyle{plain}
\newtheorem{theorem}{Theorem}[section]
\newtheorem{proposition}[theorem]{Proposition}
\newtheorem{lemma}[theorem]{Lemma}
\newtheorem{corollary}[theorem]{Corollary}

\newtheorem{hypothesis}[theorem]{Main hypothesis}
\newtheorem{conjecture}[theorem]{Conjecture}
\theoremstyle{definition}
\newtheorem{definition}[theorem]{Definition}
\newtheorem{remark}[theorem]{Remark}
\newtheorem{example}[theorem]{Example}
\endgroup

\def \div {\mathop {\rm div}\nolimits}

\def \dist {\mathop {\rm dist}\nolimits}

\def \de {\mathrm{d}}

\def \eps {\varepsilon}
\def \Om {\Omega}

\def \la {\lambda}
\def \cp {\mathrm{cap}}

\title{On a Cheeger--Kohler-Jobin inequality}
\author[I. Lucardesi]{Ilaria Lucardesi}
\address[Ilaria Lucardesi]{Dipartimento di Matematica e Informatica ``U. Dini'',
	Universit\`a di Firenze,
	Viale Morgagni 67/A, 50134 Firenze, Italy}
\email{ilaria.lucardesi@unifi.it}
\author[D. Mazzoleni]{Dario Mazzoleni}
\address[Dario Mazzoleni]{Dipartimento di Matematica ``F.~Casorati'', Universit\`a di Pavia,  Via Ferrata 5, 27100 Pavia, Italy}
\email[Dario Mazzoleni]{dario.mazzoleni@unipv.it}
\author[B. Ruffini]{Berardo Ruffini}
\address[Berardo Ruffini]{Dipartimento di Matematica\\
	Universit\`a di Bologna\\
	Piazza di Porta San Donato 5, 40126 Bologna, Italy }
\email[Berardo Ruffini]{berardo.ruffini@unibo.it}

\begin{document}
	
	\begin{abstract}
		We discuss the minimization of a Kohler-Jobin type scale-invariant functional among open, convex, bounded sets, namely
		\[
		\min\Big\{ T_2(\Omega) ^{\frac{1}{N+2}}h_1(\Omega) : \Omega\subset\R^N,\text{ open, convex, bounded}\Big\}\,
		\]
		where $T_2(\Omega)$ denotes the torsional rigidity of a set $\Omega$ and $h_1(\Omega)$ its Cheeger constant. 
		We prove the existence of an optimal set and we conjecture that the ball is the unique minimizer.
		We provide a sufficient condition for the validity of the conjecture, and an application of the conjecture to prove a quantitative inequality for the Cheeger constant. We also show lack of existence for the problem above among  several other classes of sets. As a side result we discuss the equivalence of the several definitions of Cheeger constants present in the literature and show a quite general class of sets for which those are equivalent.  
	\end{abstract}
	
	\thanks{{\bf Acknowledgments.} 
		The authors are grateful to L. Brasco and E. Parini for useful discussions on the topic of the paper. 
		The authors are members of INdAM-GNAMPA. D.~M. and B.~R. have been partially supported by the INdAM-GNAMPA project ``PACE''.  I.~L. has been partially supported by the INdAM-GNAMPA project 2022 ``Analisi variazionale via blow-up di funzionali di curvatura''. D.~M. has been partially supported by the MIUR-PRIN 2020 Mathematics for Industry 4.0.
	}

	\maketitle

	{\small
		
		\bigskip
		\noindent\keywords{\textbf{Keywords:} Cheeger constant, Kohler-Jobin inequality, quantitative estimates, Poincar\'e-Sobolev constants }
		
		\bigskip
		\noindent\subjclass{\textbf{MSC 2010:} 35P10, 39B62,
			49Q10, 
			49R05 
			}}
	\bigskip
	\bigskip

	\section{Introduction}
	Let $\Omega$ be an open, bounded set in $\R^N$. For   $1\leq r<N$ and $1\le q < Nr/(N-r)$, or for $r\ge N$ and $1\leq q<+\infty$, we define
	\begin{equation}\label{def-lpq}
		\la_{r,q}(\Om):=\inf\left\{\frac{\int_{\Omega}|\nabla u|^r\de x}{\left(\int_{\Omega}|u|^q\de x\right)^\frac{r}{q}}: u\in W^{1,r}_0(\Omega)\setminus\{0\}\right\}=\inf\left\{\frac{\int_{\Omega}|\nabla u|^r\de x}{\left(\int_{\Omega}|u|^q \de x\right)^\frac{r}{q}}: u\in C^\infty_c(\Omega)\setminus\{0\}\right\}\,,
	\end{equation}
	which can be interpreted as the principal frequency for the nonlinear eigenvalue problem \[
	-\Delta_r u=\lambda \|u\|^{r-q}_{L^q(\Om)}|u|^{q-2}u\qquad \mbox{in }\Om,\qquad\qquad u=0\qquad \mbox{on }\partial \Om\,,
	\]
	where $\Delta_r u=\div(|\nabla u|^{r-2}\nabla u)$ denotes the $r$-Laplacian of $u$.
	Alternatively, $\lambda_{r,q}(\Omega)$ can be defined as the optimal constant $C$ for the  Poincar\'e-Sobolev inequality
	\[
	C\left(\int_{\Omega} |u|^q \de x\right)^{\frac rq}\le\int_{\Omega}|\nabla u|^r\de x\,,
	\]
	in the Sobolev space  $W^{1,r}_0(\Omega)$, that is the closure of $C^\infty_c(\Omega)$ under the (semi)norm
	\[
	u\mapsto \left( \int_\Om|\nabla u|^r\de x\right)^{1/r}.
	\]
	Some of the functionals defined above have been intensively studied in the literature; one of them is the $p$-{\em torsional rigidity} (or just torsion):
	\begin{equation}\label{Tp}
	T_p(\Om):=\lambda_{p,1}^{-1}(\Om)=\left(\frac{p}{(p-1)}\max{\left\{\int_{\Om}u\, \de x - \frac1p \int_{\Om}|\nabla u|^p\de x\;:\;u\in W^{1,p}_0(\Om)\right\}}\right)^{p-1}\,.
	\end{equation}
	
	The equivalence between this definition and that in \eqref{def-lpq} can be shown by homogeneity.
	On the other hand, $\lambda_{2,2}(\cdot)$ is the classical first eigenvalue of the Dirichlet-Laplacian and, if $\Omega$ has regular enough boundary, $\lambda_{1,1}(\cdot)$ coincides with the {\it Cheeger constant} of $\Omega$,  
	\[
	h_1(\Omega):=\inf\left\{\frac{P(E)}{|E|}: E\subset\overline\Omega \right\}\,,
	\]
	where $P(\cdot)$ denotes the  Caccioppoli-De Giorgi's perimeter, and $|\cdot|$ is the $N$-dimensional Lebesgue measure (see Section~\ref{section2}).
	The P\'olya-Szeg\"o inequality, or the isoperimetric inequality if $r=q=1$, ensures that balls minimize $\lambda_{r,q}(\cdot)$ among sets of prescribed measure, that is (in a scale invariant form)
	\begin{equation}\label{quantitative}
		\lambda_{r,q}(\Omega)|\Om|^{\frac{r}{N}+\frac{r}{q}-1}\ge\lambda_{r,q}(B)|B|^{\frac{r}{N}+\frac{r}{q}-1}\,,
	\end{equation}
	for any ball $B$. Moreover, equality holds if and only if $\Omega$ is a ball, up to a set of null $r-$capacity.
	Of course, this entails that balls maximize the $p-$torsional rigidity under measure constraint, that is (in a scale invariant form)
	\[
	T_p(\Omega)|\Om|^{-\frac{p+N(p-1)}{N}}\le T_p(B)|B|^{-\frac{p+N(p-1)}{N}}\,,
	\]
	for any ball $B$. The latter is known as \emph{Saint-Venant inequality}.
	
	In this paper we will deal mostly with the functionals $T_2$ and $\lambda_{1,1}$ (which is equal to the Cheeger constant $h_1$ in the class of sets that we consider).  We stress that many results hold for more general cases. 
	
	\vspace{.5cm}	
	
	P\'olya and Szeg\"o in \cite{ps} conjectured that the product of the torsional rigidity (raised to a suitable power) and the first eigenvalue of the Dirichlet-Laplacian was minimized by balls. Intuitively, this tells that the minimality of balls for the eigenvalue is somehow more stable compared to their maximality for the torsion.
	The conjecture was proved to be true by Kohler-Jobin, who showed, in \cite{KJ,KJ1},  that
	\begin{equation}\label{original-KJ}
		T_2(\Omega)^{\frac{2}{N+2}}\lambda_{2,2}(\Omega)\ge T_2(B)^{\frac{2}{N+2}} \lambda_{2,2}(B)\,,\qquad\text{for all open, bounded }\Omega\subset \R^N\,,
	\end{equation}
	with equality if and only if $\Omega$ is a ball, up to a set of $2$-capacity zero. The exponent $2/(N+2)$ is chosen so that the functional is scale invariant.
	The  P\'olya and Szeg\"o's conjecture can   be naturally extended to the more general family of functionals $T_p^\theta(\cdot)\lambda_{p,q}(\cdot)$. This nonlinear version of it was proved by Brasco in \cite{Bra}, where he shows, with somehow simplified arguments, that 
	\[
	T_p^{\theta}(\Omega) \lambda_{p,q}(\Omega)\ge T_p^{\theta}(B) \lambda_{p,q}(B)\,,\qquad\text{for all open, bounded }\Omega\subset \R^N\,,
	\]
	whenever $B$ is a ball, $1<p<+\infty$, and $1< q < Np/(N-p)$ if $p<N$, $1< q<+\infty$ if $p\ge N$. Again, $\theta=\theta(N,p,q)$ is chosen to make the functional homogeneous and equality holds only if $\Omega$ is a ball, up to a set of $p$-capacity zero.
	We notice that the above class of parameters does not include the case $T_p^\theta(\cdot)\la_{r,q}(\cdot)$ with $p\not=r$  and in particular the interesting case $p=2,$ $r=q=1$ which involves the Cheeger constant and will be the topic of this work.
	This latter choice of parameters has a substantial difference from the other cases: for $r>1$ it is not difficult to show the existence of a minimizer in the definition \eqref{def-lpq} of $\lambda_{r,q}(\Omega)$, while for $r=1$ this is not the case, since a minimizing sequence may relax to a function with jumps, not belonging to $W^{1,1}(\Omega)$, so that a minimizer must be searched in the space of functions with bounded variation $BV(\Omega)$. Moreover, the positive level sets of the $BV$ minimizers turn out to be optimal sets for $h_1(\Omega)$ (more details on the Cheeger constant are given in Section \ref{sec:cheeger} below).

We stress here that the existence of a minimizer for the functional $T_p^{\theta}(\cdot) \lambda_{1,q}(\cdot)$ is false if stated  among merely bounded open sets of $\R^N$, due to the fact that the functional $\la_{1,q}(\cdot)$ is sensible to nontrivial $1$-capacitary modifications of a set (see Section~\ref{section2} for a survey on capacities), while $T_p(\cdot)$ to $p$-capacitary ones.
	More precisely: it is always possible to find a (compact) set $K\subset B$ with $\cp_1(K)=0$ and $\cp_p(K)>0$ for $p>1$ so that $\la_{1,q}(B\setminus K)=\la_{1,q}(B)$ for all $q\in [1,+\infty)$, but $T_p(B)>T_p(B\setminus K)$. 
	
	\noindent It is also possible to find a sequence of compacts $K_n$ with the same property above so that $T_p(B\setminus K_n)\rightarrow 0$ as $n\rightarrow \infty$, while $\la_{1,q}(B\setminus K_n)=\la_{1,q}(B)$. Therefore we have \[
	\inf{\left\{T_p(\Omega) ^{\theta}\la_{1,q}(\Omega)\;:\;\Om\subset \R^N,\;\mbox{open, bounded}\right\}}=0,
	\]
	and the infimum is not attained by any open, bounded set.
	
\noindent This same holds for all functionals $T_p(\cdot)\la_{r,q}(\cdot)$ whenever $p\not=r$, for all $q\in [1,+\infty)$.

	In this paper we aim to study the extension of \eqref{original-KJ} to the functionals $T_2^{\frac{1}{N+2}}(\cdot)\lambda_{1,1}(\cdot)$, among convex sets (and in this class, $\lambda_{1,1}$ coincides with $h_1$). 
	Restricting to the class of convex sets yields a well-posed problem, as we show in our first result below.
	
\begin{theorem}\label{thm:existence}
There exists a minimizer for the problem 
\[
\min\Big\{ T_2(\Omega) ^{\frac{1}{N+2}}h_1(\Omega) : \Omega\subset\R^N,\text{open,\ convex,\ bounded}\Big\}\,.
\]
\end{theorem}
	
	A definitely major problem is to show that  also in this case the ball is the unique minimizer. We are not able to accomplish this task and we propose it as an open problem.
	
	\begin{conjecture}\label{conjecture}[Cheeger-Kohler-Jobin inequality]
		Let $\Omega\subset \R^N$ be an open, bounded convex set and $B\subset \R^N$ be any ball. The following inequality holds true:
		\begin{equation}\label{intro-KJq}
			T_2(\Omega) ^{\frac{1}{N+2}}h_1(\Omega) \geq T_2(B)^{\frac{1}{N+2}}h_1(B)\,.
		\end{equation}	
		Moreover, equality holds in \eqref{intro-KJq} if and only if $\Omega$ is equal to a ball.
	\end{conjecture}
	
	In order to justify such a conjecture, we provide a proof of it subordinated to   a technical assumption. For the precise definition and motivation of the assumption we refer to Section \ref{sec-KJ}. Here we limit ourselves to explain its nature. 
	
	The proof of the Kohler-Jobin inequality makes use of a quite deep rearrangement technique which, after the work of Brasco in \cite{Bra}, holds as long as the parameter $p$ in the definition of $\lambda_{p,q}$ is strictly greater than $1$. Such a rearrangement is done in a way such that, given a suitably chosen function $u$:
	\begin{enumerate}
		\item the $L^p$ norm of the gradient (i.e. the numerator in the Rayleigh quotient for $\la_{p,q}$) does not increase as $u$ is rearranged;
		\item the $L^q$ norm of the function (i.e. the denominator in the Rayleigh quotient for $\la_{p,q}$) increases as $u$ is rearranged;
		\item the set $\{u>0\}$ is transformed into a ball with torsion lower than or equal to $T_2(\Omega)$.
	\end{enumerate}
	By suitably choosing the function $u$, this leads immediately to the Kohler-Jobin (and Brasco's) inequality. 
	
	In this paper via an approximation argument we are able to construct an abstract rearrangement satisfying $(1)$ and $(2)$, but we are not able  to show that $(3)$ is   satisfied too (though we believe so).
	
	The very precise statement of the hypothesis is given in \eqref{eq:hp}, in Section \ref{sec-KJ}.
	
	\begin{remark}
		Notice that the nature of the Cheeger constant may suggest to use level-by-level rearrangements, for instance by exploiting coarea type formulas. Nevertheless such an approach appears   doomed to fail, as the Kohler-Jobin rearrangement changes the height of each rearranged level set of a function. 
	\end{remark}
	
	The subordinated result we have is the following.
	\begin{theorem}\label{thm:ckjunderassumption}
		Let $\Omega\subset \R^N$ be an open, bounded convex set.
		If there exists a rearrangement $u\mapsto u^\sharp$ satisfying Hypothesis \ref{mainhp}, 
		then the minimality of the ball stated in Conjecture~\ref{conjecture} holds.
	\end{theorem}
\begin{remark}
	Note that we are not able to show, even subordinated to Hypothesis \ref{mainhp}, anything about the uniqueness of  minimizers.
\end{remark}

As a consequence of Conjecture~\ref{conjecture}, we are able to offer a new proof of the following quantitative version of \eqref{quantitative} for $r=q=1$.
	\begin{corollary}\label{thm1}
		Subordinated to the validity of  Hypothesis \ref{mainhp}, there exists a dimensional constant $\sigma(N)>0$ such that for any $\Om\subset\R^N$ open, bounded convex set and for any ball $B\subset\R^N$, we have
		\begin{equation}\label{qscheeger}
			|\Omega|^{\frac{1}{N}} h_1(\Omega) - |B|^{\frac{1}{N}} h_1(B)\geq \sigma (N)\alpha(\Omega)^2\,,
		\end{equation}
		where $\alpha(\Om)$ denotes the  Fraenkel asymmetry of $\Om$, see~\eqref{fraenky}.
		Moreover, the power $2$ in~\eqref{qscheeger} is sharp, in the sense that it can not be replaced by any lower number. 
	\end{corollary}
	This improvement of the Cheeger inequality was first showed by Figalli, Maggi, and Pratelli, among open sets: in~\cite{fimp}  they provide a short proof of this fact, based on the quantitative version of the isoperimetric inequality \cite{fmp,fimp2,cl}.
	
	The  approach in this paper is quite different, and borrows an idea from  \cite{bdpv}, where Brasco, De Philippis, and Velichkov showed that $\lambda_{2,2}(\cdot)$ satisfies a (asimptotically sharp) quantitative estimate of the form
	\[
	|\Omega|^{\frac{2}{N}}\lambda_{2,2}(\Omega)-|B|^{\frac{2}{N}}\lambda_{2,2}(B)\ge C(N)\alpha(\Omega)^2\,,
	\]
	by relating the stability of $\lambda_{2,2}(\cdot)$ to that of the  torsional rigidity $T_2(\cdot)$.
	In the very same spirit, we obtain~\eqref{qscheeger} by combining the Cheeger--Kohler-Jobin inequality \eqref{intro-KJq}, together with the quantitative stability of the $2$-torsional rigidity provided in \cite{bdpv}.
	Though this result is weaker than the one in~\cite{fimp}, we believe it is an interesting application of Conjecture~\ref{conjecture}, and it should be a stimulus toward proving (or disproving) it.
	
	{\bf Plan of the paper. }
	The article is organized as follows. After some preliminaries in Section~\ref{section2} about the geometric measure theory notions needed to define the Cheeger constant and some recalls on capacities, in Section~\ref{sec:cheeger} we  survey the (several) different   definitions of the Cheeger constant appearing in literature. This is in some sense of independent interest, even if it does not have any claim of novelty.
	Section~\ref{sec:proofthmexistence} is devoted to the proof of Theorem~\ref{thm:existence}. After presenting in Section~\ref{sec-KJ} the Kohler-Jobin rearrangement and the main Hypothesis \ref{mainhp}, in Section~\ref{section6} we discuss Conjecture~\ref{conjecture} and then we prove Theorem~\ref{thm:ckjunderassumption}.
	Section~\ref{sec-qe} deals with the proof of the quantitative estimate of Corollary~\ref{thm1}.

	\section{Preliminaries}\label{section2}
	
	In this section we collect some well-known facts on geometric measure theory, which will serve our scopes later in the paper. We refer to~\cite[Chapter~3]{afp} for more details.
	\subsection{The perimeter and its properties}
	The measure theoretic perimeter (shortly: perimeter) of a Borel set $E\subset\R^N$ is the quantity
	\[
	P(E):=\sup\left\{\int_{E}\nabla\cdot\phi \, \de x \,:\,\phi\in C^1_c(\R^N,\R^N), \,\, \|\phi\|_\infty\le 1 \right\}\,.
	\]
	If $P(E)<+\infty$ then we say that $E$ has finite perimeter. 
	Equivalently, it can be defined in the setting of  functions of bounded variation as the total variation of the distributional derivative of characteristic functions. We recall that if $\Om\subset\R^N$ is an open set, $u\in L^1(\Om)$ is a function of bounded variation, and we write $u\in BV(\Om)$, when the distributional derivative $Du$ of $u$ is an $\R^N$-valued finite Radon measure.
	If $E\subset\R^N$ is a set of finite perimeter, then $\chi_E\in BV(\R^N)$ and 
	\begin{equation}\label{eq:perdistrib}
	P(E)=|D\chi_E|(\R^N)=:\|D\chi_E\|_{TV(\R^N)}.
	\end{equation}	
	Whenever it exists, the quantity
	\[
	[0,1]\ni\theta_E(x):=\lim_{r\to0}\frac{|E\cap B_r(x)|}{|B_r(x)|}\,,
	\]
	is called the density of a Borel set $E$ at $x$. We denote by $E^t$ the subset of points of $\R^N$ such that $\theta_E(x)=t$, and we call essential boundary of $E$ the set $\partial^e E=E\setminus (E^0\cup E^1)$. Eventually, we define the reduced boundary of $E$ as the set $\partial^*E\subset\partial^eE$ of points of the essential boundary such that the measure theoretic inner unit normal 
	\[
	\nu_E(x):=\lim_{r\to0}\frac{D\chi_E(B_r(x))}{|D\chi_E|(B_r(x))}
	\] 
	exists.
	
	The geometry  of the boundary of sets of finite perimeter is described in the two cornerstones of geometric measure theory: the De Giorgi's and the Federer's structure theorems.
	\begin{theorem}[De Giorgi's Structure Theorem]
		Let $E$ be a set of finite perimeter. Then $\partial^*E$ is $\mathcal H^{N-1}-$rectifiable and $P(E)=\mathcal H^{N-1}(\partial^*E)$. Moreover, if $x\in \partial^*E$, then $(E-x)/r$ converges in $L^1_{\rm loc}$ to the hyperspace orthogonal to $\nu_E(x)$, as $r\to 0$. Eventually, the following divergence formula holds true
		\[
		\int_{E}\nabla\cdot\phi\, \de x=-\int_{\partial^*E}\phi\cdot\nu_E \, \de \mathcal H^{N-1}(x)\,,
		\]  
		for any vector field $\phi\in C^1_c(\R^N,\R^N)$.
	\end{theorem}
	
	\begin{theorem}[Federer's Structure Theorem]
		Let $E$ be a set of finite perimeter. Then $\partial^*E\subset E^{1/2}$ and $\mathcal H^{N-1}(\partial^eE\setminus\partial^*E)=0$. In particular $\partial^*E$, $E^{1/2}$, and $\partial^eE$ are equivalent, up to a $\mathcal H^{N-1}-$negligible set. 
	\end{theorem}

	\subsection{The isoperimetric inequality and its quantitative improvement}
	
	The Kohler-Jobin inequality is based on the simple principle of slicing  the energy functionals defining $\lambda_{p,q}$ and $T_p$ horizontally, then rearranging the level sets of the involved functions in a suitable way.
	The energies before and after the rearrangement are compared by exploiting the isoperimetric inequality: for any set $E$ of finite measure, 
	\[
	P(E)-P(B)\ge0
	\]
	whenever $B$ is a ball of measure $|E|$, with equality if and only if $E$ coincides with $B$ up to a negligible set. Equivalently, by homogeneity, for all $E\subset\R^N$ of finite measure it holds
	\[
	P(E)- N\omega_N^{1/N}|E|^{(N-1)/{N}}\ge0\,,
	\]
	where $\omega_N$ denotes the measure of the ball with unit radius in $\R^N$, and again equality holds if and only if $E$ is a ball.
A stronger version of the isoperimetric inequality actually holds, see for example~\cite{fmp}: there exists a dimensional constant $C_N$ such that, for any set $E\subset\R^N$ of finite measure, it holds
	\begin{equation}\label{isoquantitativa}
		\frac{P(E)-N\omega_N^{1/N}|E|^{(N-1)/{N}}}{N\omega_N^{1/N}|E|^{(N-1)/{N}}}\ge C_N \alpha(E)^2,
	\end{equation}
	where  $\alpha(E)$ is the \emph{Fraenkel asymmetry} of the set $E$, namely 
	\begin{equation}\label{fraenky}
		\alpha(E):=\inf_{x\in\R^N}\left\{\frac{|E\Delta (B+x)|}{|\Omega|}: \text{ $B\subset\R^N$ is a ball,  }|B|=|E|  \right\},
	\end{equation}
and $U\Delta V$ stands for the symmetric difference between the sets $U$ and $V$.\\
	It is worth stressing that the exponent $2$ in the quantitative estimate \eqref{isoquantitativa} is \emph{sharp}, in the sense that it can not be replaced by any lower number. 
	
	\subsection{Sobolev capacities}
	We   briefly recall the definition of the $p$-capacity of a measurable set $E\subset \R^N$, for $1\leq p<N$. We refer for more properties and details about it to~\cite{eg}.
	We define (here $A^\circ$ stands for the open interior of a set $A\subset\R^N$)
	 \[
	\cp_p(E):=\inf\left\{\int_{\R^N}|\nabla u|^p\,dx\;:\;u\geq 0,\;u\in L^{p^*}(\R^N),\; \nabla u\in L^p(\R^N;\R^N),\;\{u\geq 1\}^\circ \supset E\right\},
	\]
	and clearly for all $K\subset\R^N$ compact it becomes \[
	\cp_p(K)=\inf\left\{\int_{\R^N}|\nabla u|^p\,dx\;:\;u\in C^\infty_c(\R^N),\;u\geq \chi_K\right\}.
	\]
	The only properties that we highlight here are the link of the $p$-capacity with Hausdorff and Lebesgue measures of a set.
	Namely that there exist constants $C_1(N,p),C_2(N,p)>0$ such that, for all $E\subset\R^N$\[
	\cp_p(E)\leq C_1 \mathcal H^{N-p}(E),\qquad \mathcal L^{N}(E)\leq C_2 \cp_p(E)^{\frac{N}{N-p}},\qquad \text{ for }1\leq p<N.
	\]
	
	\section{On the different definitions of the Cheeger constant}\label{sec:cheeger}
	A main role in this paper is played by the Cheeger constant of a set. Such a constant was introduced in~\cite{cheeger} to obtain lower bounds for the first eigenvalue of the Laplace-Beltrami operator. While its original definition is given on Riemannian compact manifolds without boundary, it has lately found many applications in the Euclidean setting, where it can be defined in several ways. Here we briefly survey such different definitions (for more details, see \cite{surveyleonardi}), and we offer a criterion under which all the corresponding constants coincide. This allows us to switch from one definition to the other in the rest of the paper. 
	\begin{definition}
		Let $\Omega$ be an open, bounded set in $\R^N$. Then the Cheeger constant is either
		\begin{itemize}
			\item[]
			\[
			h_1(\Omega):=\inf\left\{\frac{P(E)}{|E|}:E\subset\overline{\Omega}\right\},\text{ or }
			\]
			\item[]
			\[
			h(\Omega):=\inf\left\{\frac{P(E)}{|E|}:E\Subset {\Omega}\right\},\text{ or }
			\]
			\item[]
			\[
			\lambda_1(\Omega):=\inf\left\{\frac{\|Du\|_{TV(\R^N)}}{\|u\|_{L^1(\Omega)}}\,:\,u\in BV(\overline\Omega)\setminus \{0\},\,u|_{{\R^N\setminus\overline{\Omega}}}=0\right\},\text{ or }
			\]
			\item[]
			\[
			\lambda_{1,1}(\Omega):=\inf\left\{\frac{\|\nabla u\|_{L^1(\R^N)}}{\|u\|_{L^1(\Omega)}}:u\in C^{\infty}_c(\Omega)\setminus \{0\}\right\}\,.
			\]
		\end{itemize}
		We say that $E_\Om$ is a \emph{Cheeger set} for $\Om$ if $h_1(\Om)=\frac{P(E_\Om)}{|E_\Om|}$ and that $\Om$ is \emph{self--Cheeger} if $\Om=E_\Om$ up to sets of $1$-capacity zero.
	\end{definition}
The functional $h_1$ is the one in the spirit of the original work by Cheeger~\cite{cheeger}.  
Now, we can see that two of the above infima are actually minima.
\begin{lemma}\label{le:existcheeger}
			Let $\Omega$ be an open, bounded set of $\R^N$, then the infima in the definition of $h_1$ and $\la_1$ are actually minima. Moreover, $\lambda_1(\Om)$ is attained by a characteristic function.
		\end{lemma}
		\begin{proof}
We show the existence of a minimum for $h_1$, the case of $\la_1$ is similar, see~\cite[Theorem~8]{frikaw} for the last claim. 
It is clear that $h_1(\Omega)<+\infty$. 
			Let $(E_n)_n\subset \overline\Omega$ be a minimizing sequence for $h_1(\Omega)$ such that \[
			 \frac{P(E_n)}{|E_n|}\leq h_1(\Omega)+\frac1n.
			\] 
Then, using the isoperimetric inequality, we obtain\[
			\begin{split}
				N\omega_N^{1/N}|E_n|^{1-\frac1N}\leq P(E_n)\leq (h_1(\Omega)+\frac1n)|E_n|,\qquad \text{hence,}\qquad |E_n|^{1/N}\geq C.
			\end{split}
			\]
On the other hand, \[
P(E_n)\leq (h_1(\Om)+\frac1n)|E_n|\leq (h_1(\Om)+\frac1n)|\Omega|,\qquad \text{hence,}\qquad P(E_n)\leq C.
\]
			By the equiboundedness of the perimeter, we can find a nonrelabeled subsequence and a nonempty set of finite perimeter $E\subset \overline \Omega$ such that $E_n\to E$ in $L^1$. By the lower semicontinuity of the perimeter, we conclude \[
\frac{P(E)}{|E|}\leq 			\liminf_{n\to+\infty}\frac{P(E_n)}{|E_n|}= h_1(\Omega)\,.\qedhere
			\] 
\end{proof}

By coarea formula, it is not difficult to show that any minimizer in the definition of $\la_1$ has the property that each of its level sets is a Cheeger set of $\Omega$, that is, a minimizer of $h_1(\Omega)$, see \cite{frikaw}. 
	Nonetheless, the above defined constants are not the same in general. 
	{\begin{example}
			Let us consider $B_1$, the ball centered at $0$ of radius $1$ in $\R^2$. Let $U=[-\tau,\tau]\times\{0\}$ for some $\tau\in(0,1)$ and $\Omega=B_1\setminus U$. 
			One can easily show by the minimality of the ball for the Cheeger constant (see \eqref{quantitative}) that $h(\Omega)>h(B_1)$. By a more delicate argument one might show also that 
			\[
			h(\Omega)\ge \frac{P(B_1)}{|B_1|}+2\mathcal H^1(U)\ge h_1(B_1)+2\tau\,
			\]
			as long as $\tau$ is small enough.
			The very same example holds for $\lambda_{1,1}$ in place of $h$.
		\end{example}
	} 
		Notice that in the previous examples  we removed a {\em (N-1)-dimensional   manifold}  from a regular set (the ball). We wish to investigate if this condition is somehow close to be sharp. The answer happens to be positive, as shown in the next proposition. 
	
	\begin{proposition}
		Let $\Omega$ be an open, bounded set of $\R^N$ such that
		\[
		P(\Omega)=\mathcal H^{N-1}(\partial\Omega)<+\infty\,.
		\]
		Then $h_1(\Omega)=h(\Omega)=\lambda_1(\Omega)=\lambda_{1,1}(\Omega)$.
	\end{proposition}
	\begin{proof} We claim that, since $P(\Omega)=\mathcal H^{N-1}(\partial \Omega)$, the same happens for a minimizer $E$ of $h_1(\Omega)$, which exists by Lemma~\ref{le:existcheeger}. To show this, we first split the reduced boundary of $E$ in the following way
		\[
		\partial^*E=(\partial^*E\cap\partial\Omega) \cup (\partial^*E\setminus\partial \Omega)\,.
		\]
		From De Giorgi's Theorem, it is well known that outside the contact points, i.e. in $\partial^*E\setminus\partial\Omega$, the set $E$ is regular, so that $\mathcal H^{N-1}(\partial^*E\setminus\partial\Omega)=\mathcal H^{N-1}(\partial E\setminus\partial\Omega)$, see for example~\cite[Proposition~3.5~(iv)]{surveyleonardi}. 
		As for the contact points, we clearly have $\mathcal H^{N-1}(\partial^*E\cap \partial\Omega)\le\mathcal H^{N-1}(\partial E\cap \partial \Omega)$.
		On the other hand, $\mathcal H^{N-1}-$a.e. $x\in \partial E\cap \partial\Omega$ belongs to $\partial E\cap \partial^*\Omega$, from our hypothesis on $\Omega$ and Federer's Theorem.
		Finally, thanks to \cite[Prop. 3.5, point (vii)]{surveyleonardi}, any $x\in \partial^*\Omega\cap\partial E$ belongs to $\partial^*E$ (namely: $\partial E$ meets $\partial \Omega$ tangentially), so that $\mathcal H^{N-1}(\partial^* E\cap \partial\Omega)=\mathcal H^{N-1}(\partial E\cap \partial \Omega)$.  
		Summing up all the informations, we have
		\[
		\begin{aligned}
			P(E)&=\mathcal H^{N-1}(\partial^*E)\\
			&=\mathcal H^{N-1}(\partial^*E\cap\partial \Omega)+\mathcal H^{N-1}(\partial^*E\setminus\partial \Omega) =\mathcal H^{N-1}(\partial^*E\cap\partial \Omega)+\mathcal H^{N-1}(\partial E\setminus\partial \Omega)\\
			&=\mathcal H^{N-1}(\partial E\cap\partial \Omega)+\mathcal H^{N-1}(\partial E\setminus\partial \Omega)=\mathcal H^{N-1}(\partial E)\,.
		\end{aligned}
		\]
		This shows the  claim.
		From this fact, we can exploit~\cite[Theorem 1.1]{schmidt}, which ensures the existence of a sequence of smooth sets $E_n$ compactly contained inside $E$, which approximate $E$ both in $L^1$ and in perimeter. In particular,
		\[
		h_1(\Omega)\le  h(\Omega)\leq \lim_{n\to+\infty}\frac{P(E_n)}{|E_n|}=\frac{P(E)}{|E|}=h_1(\Omega),
		\]
		namely $h_1(\Om)=h(\Om)$.
	
Taking the same approximating sequence above $E_n\Subset E\subset \Omega$,  we construct a sequence $(u_n)_n \subset C^\infty_c(\Omega)$ such that 
\begin{equation}\label{eq:approx}
		\frac{\int_{\Omega}|\nabla u_n|\, \de x}{\int_{\Omega}|u_n|\, \de x}=\frac{P(E)}{|E|}+o_n(1).
\end{equation}
The sequence $(u_n)_n$ is defined as follows: we take $u_n:=\rho_{\varepsilon/2}*(\chi_{E_n^\varepsilon})$, where $A^\varepsilon$ is the set of points of $A$ whose distance from $\partial A$ is larger than $\varepsilon$, and $\rho_t$ is a positive mollifying kernel of total mass $1$ and such that the support of $\rho_1$ is contained in a ball of radius $1$. Notice that such a construction is admissible since $\mathrm{dist}(\partial E,\partial E_n)>0$. Moreover, the sequence satisfies \eqref{eq:approx}.
Recalling the definitions of $\la_1$ and $\la_{1,1}$,~\eqref{eq:approx} implies that $\lambda_1(\Omega)\leq \lambda_{1,1}(\Omega)\leq h(\Omega)=h_1(\Omega)$.

In order to conclude the proof, it is then enough to show that $\lambda_{1}(\Omega) \geq h_1(\Omega)$.
To this aim, we recall that, thanks to Lemma~\ref{le:existcheeger}, $\la_1(\Omega)$ is attained by a characteristic function of a set $F\subset\overline  \Omega$, hence, recalling also the characterization of the perimeter~\eqref{eq:perdistrib},\[
\la_1(\Omega)=\frac{\|D\chi_F\|_{TV(\R^N)}}{\|\chi_F\|_{L^1(\Omega)}}=\frac{P(F)}{|F|}\geq h_1(\Omega)\,.\qedhere
\] 
	\end{proof}
We note that if $\Omega\subset \R^N$ is open, bounded and convex, then $P(\Omega)=\mathcal H^{N-1}(\partial \Omega)$, therefore all the definitions of Cheeger constant are equivalent. In the rest of the paper, since we will be dealing with convex sets, we are allowed to use any of the definitions, depending on what is more convenient.
	\begin{remark}
		Notice that, for sets of finite perimeter, the inequality $P(E)\le\mathcal H^{N-1}(\partial E)$ is always true, but the equality does not hold in general as long as the $\mathcal H^{N-1}-$ measure of $E^0\cap\partial E$ and $E^1\cap\partial E$ is nonzero, as a consequence of Federer's Theorem. 
		The condition $\mathcal H^{N-1}(E^0\cap \partial E)>0$ is quite pathological. Indeed, due to the fact that sets of finite relative perimeter satisfy density estimates on their boundary, whenever this happens to be true, then $E$ can not even support a relative isoperimetric inequality. For a proof see \cite[Lemma 3.5]{saracco}. On the other hand, the condition $\mathcal H^{N-1}(E^1\cap\partial E)>0$ can hold even for self--Cheeger sets, as shown in~\cite[Section~2]{leonardisaracco}.
	\end{remark}
	
	Some important and difficult topics related to the Cheeger problem are, given an open, bounded set $\Om\subset\R^N$, to study if its Cheeger set is unique, its regularity and whether $\Om$ is self--Cheeger or not.
	In particular, for our purposes, it is important to know that if $\Om$ is convex, then $E_\Om$ is unique and has the same regularity of $\Om$. This result  has been proved by Alter and Caselles~\cite[Theorem~1]{ac}, we recall it for completeness.
	\begin{theorem}[Alter-Caselles]
		There is a unique Cheeger set inside any non-trivial open, bounded convex set in $\R^N$.
		The Cheeger set is convex and of class $C^{1,1}$.
	\end{theorem}

	\section{Proof of Theorem~\ref{thm:existence}}\label{sec:proofthmexistence}
	
	This section is devoted to the proof of the existence result of a minimizer for the shape functional $T_2^{\frac{1}{N+2}}h_1$ among bounded open convex sets.

	A key tool will be the adaptation of the well-known John's Lemma \cite{FJ}. We state it on the class of open, bounded convex sets, though it is more often stated in the class of convex bodies. Clearly the two formulations are equivalent.
	
	Here and in the rest of the paper we say that $A\lesssim B$ if there exists a constant $C$ depending only on the dimension $N$ such that $A\le C B$. We say that $A\simeq B$ if $A\lesssim B\lesssim A$.
	\begin{lemma}\label{John}
		Let $\Omega\subset \mathbb R^N$ be an open, bounded   convex set. Then there exists a (open) parallelepiped $K$ such that
		$$
		K \subset \Omega \subset C(N) K,
		$$ where $C(N)>0$ is a constant which only depends on the dimension $N$. In particular, $|\Omega|\simeq |K|$ and $\mathrm{diam}(\Omega)\simeq \mathrm{diam}(K)$.
	\end{lemma}

	\begin{proof}[Proof of Theorem \ref{thm:existence}]
		Let $(\Omega_n)_n$ be a minimizing sequence for $T_2^{\frac{1}{N+2}}h_1$ among open convex sets of $\mathbb R^N$.
		Since the functional is scale invariant, we may assume that $|\Omega_n|=1$ for every $n\in \mathbb N$. 
		Let us prove that the diameters $\left(\mathrm{diam}(\Omega_n)\right)_n$ are uniformly bounded. Assume by contradiction that this is not true, hence there exists a subsequence (not relabeled) with diameters diverging to $+\infty$. 
		
		Let us consider the associated sequence $(K_n)_n$ of parallelepipeds provided by Lemma \ref{John}. By construction, we have that
		$$
		|K_n|\simeq 1, \quad \mathrm{diam}(K_n) \to +\infty.
		$$
		It is now convenient to rewrite $K_n$ as a scaling of a parallelepid of volume 1, as follows:
		$$
		K_n=|K_n|^{1/N} Q_n \quad \hbox{with}\quad |Q_n|=1.
		$$
		Without loss of generality (up to a rotation and a translation), we may take $Q_n$ of the form
		$$
		Q_n:=(0,\ell_1^{(n)})\times \ldots \times (0, \ell_N^{(n)}),
		$$
		with 
		$$
		\Pi_{i=1}^N \ell_i^{(n)} = 1
		$$
		and with ordered sides
		$$
		\quad 0 < \ell_i^{(n)}\leq \ell_j^{(n)}\qquad \forall \;1\leq i<j\leq N.$$
		The assumption $\mathrm{diam}(K_n) \to +\infty$ entails that also $\mathrm{diam}(Q_n) \to +\infty$. This fact, together with $|Q_n|=1$, implies that
		for $n\to \infty$
		$$
		\ell_N^{(n)} \to +\infty \quad \hbox{and}\quad  
		\ell_1^{(n)}\to 0.
		$$
		
		Using the monotonicity of $T_2$ and $h_1$ with respect to set inclusion, we infer that
		\begin{equation}\label{assurdo}
			T_2^{\frac{1}{N+2}}(\Omega_n)h_1(\Omega_n)\geq \frac{1}{C(N)}T_2^{\frac{1}{N+2}}(K_n)h_1( K_n) = \frac{1}{C(N)}T_2^{\frac{1}{N+2}}(Q_n)h_1(Q_n),
		\end{equation}
		where $C(N)$ is the constant appearing in Lemma \ref{John}. 
		Let us bound from below the two functionals $h_1$ and $T_2$, computed at $Q_n$. On the one hand, we exploit the following estimate, proved in  \cite[Corollary 5.2]{brascopini}: for every $\Omega\subset \mathbb R^N$ open bounded convex set, there holds
		$$
		h_1(\Omega) \geq \frac{1}{N} \frac{P(\Omega)}{|\Omega|}.
		$$
		Taking $\Omega=Q_n$, recalling the assumption $|Q_n|=\Pi_{i} \ell^{(n)}_{i}=1$, and computing 
		$$
		P(Q_n)=\sum_{i=1}^N \Pi_{i\neq j} \ell_j^{(n)} = \sum_{i=1}^N \frac{1}{\ell_{i}^{(n)}} \geq  \frac{1}{\ell_1^{(n)}},
		$$ 
		we infer that
		\begin{equation}\label{h1n}
			h_1(Q_n) \geq \frac{1}{N} \frac{1}{\ell_1^{(n)}}.
		\end{equation}
		Let us now pass to the torsional rigidity. By definition, for every $u\in W^{1,2}_0(Q_n)\setminus \{0\}$, we have
		$$
		T_2(Q_n)\geq \frac{\left( \int_{Q_n} u(x)\, \mathrm{d} x \right)^2}{\int_{Q_n} |\nabla u(x)|^2 \, \mathrm{d}x}.
		$$
		Writing $\mathbb R^N\ni x = (x_1, \ldots, x_N)$, we take $u(x):= \Pi_{i=1}^N \sin (\pi x_i/ \ell_{i}^{(n)})$. The function $u$ clearly belongs to $W^{1,2}_0(Q_n)\setminus \{0\}$. The numerator of the Rayleigh quotient above is a constant, which only depends on $N$:
		$$
		\int_{Q_n} u(x)\, \mathrm{d}x = \Pi_{i=1}^N \int_0^{\ell_i^{(n)}}  \sin (\pi x_i/ \ell_{i}^{(n)})\, \mathrm{d}x_i = 
		\Pi_{i=1}^N \frac{\ell_i^{(n)}}{\pi} \int_0^{\pi}  \sin (t)\, \mathrm{d}t = 
		\left(\frac{2}{\pi}\right)^N \Pi_{i=1}^N \ell_i^{(n)}=  \left(\frac{2}{\pi}\right)^N.
		$$
		As for the denominator, since
		$$
		\left(\nabla u(x)\right)_i = \frac{\pi}{\ell_i^{(n)}} \cos(\pi x_i / \ell_i^{(n)}) \Pi_{j\neq i} \sin(\pi x_j/\ell_j^{(n)}), 
		$$
		we obtain
		$$
		\int_{Q_n} |\nabla u(x)|^2 \, \mathrm{d}x 
		= \pi^2 \sum_{i=1}^N \frac{1}{\left( \ell_i^{(n)}\right)^2 }\int_0^{\ell_i^{(n)}}\cos^2(\pi x_i/\ell_i^{(n)})\, \mathrm{d}x_i \Pi_{j\neq i }  \int_0^{\ell_j^{(n)}}  \sin^2(\pi x_j/ \ell_{j}^{(n)})\, \mathrm{d}x_j
		= \frac{\pi^2}{2^N}  \sum_{i=1}^N \frac{1}{\left( \ell_i^{(n)}\right)^2 }.
		$$
		All in all we obtain
		\begin{equation}\label{Tn}
			T_2(Q_n) \geq \frac{4^N}{\pi^{N+2}}\left[ \sum_{i=1}^N \frac{1}{\left( \ell_i^{(n)}\right)^2 }\right]^{-1} \geq \frac{4^N}{\pi^{N+2}}\left[ N \frac{1}{\left( \ell_1^{(n)}\right)^2 }\right]^{-1} = \frac{4^N}{N \pi^{N+2}}\left( \ell_1^{(n)}\right)^2, 
		\end{equation}
		where, in the second inequality, we have used the assumption $\ell_i^{(n)} \geq \ell_1^{(n)}$.
		
		By combining \eqref{h1n} with \eqref{Tn} we get
		$$
		T_2^{\frac{1}{N+2}}(Q_n) h_1(Q_n) \geq \tilde{C}(N)   \frac{1}{\left(\ell_1^{(n)}\right)^{\frac{N}{N+2}}}
		$$
		where for brevity we have set $\tilde{C}(N):= 4^\frac{N}{N+2}/(\pi N^{1 + \frac{1}{N+2}} )$.
		
		Letting $n\to \infty$, we conclude that $T_2^{\frac{1}{N+2}}(Q_n) h_1(Q_n)$ goes to $+\infty$. This, together with~\eqref{assurdo}, provides a contradiction. 
		
We have thus showed that the minimizing sequence $(\Omega_n)_n$ with $|\Omega_n|=1$ has uniformly bounded diameter. In order to conclude, we apply the Blaschke selection theorem: this guarantees the existence of a subsequence converging to an open, bounded convex $\Omega^*$, with respect to the complementary Hausdorff distance; with respect to such a convergence, both $T_2$ and $h_1$ are continuous (see, e.g., \cite{H} for the continuity of $T_2$ and \cite{par} for the continuity of $h_1$), implying that 
		the shape $\Omega^*$ is a minimizer. This concludes the proof.
	\end{proof}

	\section{Kohler-Jobin rearrangement technique}\label{sec-KJ}
		
	We present here a possible strategy to attack Conjecture~\ref{conjecture}, which is based on the Kohler-Jobin radial rearrangement technique, later extended by Brasco to the nonlinear case $p\ne2$. In this section we write a short explanation of this tool and its application in our setting (see Lemma \ref{apprAp}). At the same time, we introduce some notations used in the sequel. 
	
	\begin{remark}In this section we shall work with the $p$-torsion and its associated Kohler-Jobin rearrangement for general $p>1$. This will allow us to construct a peculiar rearrangement by means of a limit argument as $p\to 1$.
	\end{remark}

	The cornerstone of the Kohler-Jobin inequality is a rearrangement inequality which acts as follows: given $1<p<+\infty$ and a non-negative function $u$ in the Sobolev space $W^{1,p}_0(\Omega)$, one constructs a rearrangement $u^*$ of $u$, belonging to $W^{1,p}_0(B)$ for some ball $B$ with $T_p(\Omega)\geq T_p(B)$, such that
	\begin{equation}\label{condizioni}
		\int_{\Omega}|\nabla u|^p\de x=\int_{B}|\nabla u^*|^p\de x
		\qquad \text{ and }\qquad 
		\int_{\Omega}|u|^q\de x\le \int_B |u^*|^q\de x\,,\qquad q\in[1,+\infty)\,.
	\end{equation}
	A somehow  natural idea, in order to obtain the Kohler-Jobin inequality, is to consider the function $u^*$ such that any level set of $u^*$ is a ball  $B_{r(t)}$ centered at $0$ with
	\[
	T_p(B_{r(t)})=T_p(\{u>t\})\,.
	\]
	Unfortunately this idea fails. In particular, the second requirement in \eqref{condizioni} can not hold in general: if $\{u>t\}$ is not a ball for too many values of $t$, then
	$$
	\int_\Omega| u|^q\de x > \int_B |u^*|^q\de x\,,
	$$
	from the Saint-Venant inequality (while nothing can be said on the $L^p$-norms of the gradients).
	
	This suggests that the rearrangement must somehow take into account the other level sets of $u$. The successful idea of Kohler-Jobin was to introduce the following modification of the torsional rigidity.
	\begin{definition}[\cite{Bra,KJ}]
		Let $\Omega$ be an open, bounded set  and $1<p<+\infty$. 
		We say that $u\in W^{1,p}_0(\Omega)\cap L^\infty(\Omega)$ is a \emph{reference function for $\Om$} if $u\geq 0$ in $\Om$ and 
		\begin{equation}\label{cattiva}
			t\mapsto \frac{|\{x\in\Om\;:\;u(x)>t\}|}{\int_{\{u=t\}}|\nabla u|^{p-1}\,\de \mathcal{H}^{N-1}}\in L^\infty([0,\|u\|_{L^\infty(\Om)}])\,.
		\end{equation}
		We call $\mathcal A_p(\Om)$ the set of all ($p$-)reference functions for $\Om$.
		
		\noindent Then, for any $u\in \mathcal A_p(\Om)$, the {\em modified torsional rigidity} is the functional, depending on $p$, $\Omega$ and $u$, defined by
		\[
		T_{p,mod}(\Omega,u)=\left( \frac{p}{p-1}\sup\left\{\int_{\Omega}g\circ u\, \de x-\frac1p \int_\Omega |\nabla g\circ u|^p\de x:g\in \mathrm{Lip}[0,\|u\|_{L^\infty(\Omega)}],\,g(0)=0\right\}\right)^{p-1}.
		\] 
	\end{definition}
We note that with $g(x)=x$ in the definition of $T_{p,mod}$ we are back to the usual torsional rigidity, see~\eqref{Tp}.	
	
		The features of $T_{p,mod}$ that will be used later are collected in the next Lemma (for the proof, we refer to \cite[Proposition~3.8]{Bra} and~\cite{KJ}).
	\begin{lemma}\label{boiade}
		Let $1<p<+\infty$, $\Om\subset\R^N$ be an open, bounded set, and $u\in \mathcal A_p(\Om)$. Then
		\begin{itemize}
			\item[(i)] $T_{p,mod}(\Omega,u)\le T_p(\Omega)$;
			\item[(ii)] if $B$ is a ball such that $T_{p,mod}(\Omega,u)=T_p(B)$, then $|B|\le|\Omega|$. Equality holds in the latter if and only if $\Omega=B$ and $u$ is a radial function.
		\end{itemize}
	\end{lemma}

	The idea in \cite{KJ,KJ1,Bra} is then the following: given $u\in \mathcal A_p(\Om)$, define $u^*$ as the function such that for any $t\in[0,\|u\|_{L^\infty(\Omega)}]$, the set $\{u^*_p>\psi(t)\}$ is a ball with $p-$torsional rigidity equal to $T_{p,mod}(\{u>t\},(u-t)_+)$, where for a function $f$ we call $f_+=\max\{f,0\}$ the positive part of $f$.  Here $\psi(t)$ is a suitably chosen decreasing function. Its definition is quite implicit and it is the core of the rearrangement. Here we only stress that in general if $\Omega$ is not a ball, $\psi$ can not be the identity. With this construction, it is possible to show a rearrangement result as follows (for a proof we refer to~\cite[Proposition~4.1 and Remark~4.3]{Bra} and \cite{KJ}).

The following theorem contains the key features of the aforementioned rearrangement.
	\begin{theorem}[Kohler-Jobin Rearrangement Theorem]
		Let $1<p<+\infty$, $\Om\subset\R^N$ be an open, bounded set, $u\in \mathcal A_p(\Om)$, and $B^p$ the origin centered ball such that $T_p(B^p)=T_{p,mod}(\Omega,u)$. Then, for every $q \geq 1$, there exists a radially symmetric decreasing function $u^*_p\in W^{1,p}_0(B^p)$, such that
		\[
		\int_{\Omega} |\nabla u|^p\de x=	\int_{B^p} |\nabla u^*_p|^p\de x\quad\text{ and }\quad 	\int_{\Omega} |u|^q\de x\le\int_{B^p} |u^*_p|^q\de x\,.	
		\]
		Moreover, if $q>1$, equality holds in the latter if and only if $u$ is already a radially decreasing function.
	\end{theorem}
	In the sequel we will call such an $u^*_p$ \emph{Kohler-Jobin rearrangement of $u$}. 

		\begin{remark}
			We note that while in~\cite[Proposition~4.1]{Bra} the case $q=1$ is not included, it is observed in~\cite[Remark~4.3] {Bra} that it still holds (though the lack of strong convexity of the power does not allow to treat the equality cases).
		\end{remark}
	
	\medskip
	
	{We can not apply directly the Kohler-Jobin rearrangement for our aims, but a sort of limiting version will serve our purposes.  To this scope let us consider a function $u\in \mathcal A(\Omega)$, where
\[
\mathcal A(\Omega):= \bigcap_{r\in (1,2]}	 \mathcal A_r(\Omega).	
\]
Then, for every $r\in (1,2]$, it is well-defined 
\[
u^*_r\in W^{1,r}_0(B^r)\subset BV_0(B^r)=\{v\in BV(B^r)\,:\, \text{ $v=0$ outside of $B^r$}\,\}, \]
the Kohler-Jobin rearrangement of $u$, where $B^r$ is a ball centered at the origin such that $T_r(B^r)=T_{r,mod}(\Omega,u)$. We extend to zero $u^*_r$ outside $B^r$, if needed. 
Let us first show that the balls $B^r$ have radii uniformly bounded, thus they are all contained in a big ball $B_R$ for some $R>1$.
In fact, first of all, using the monotonicity of $T_r$ and the definition of $T_{r,mod}$, we have (using also the Saint-Venant inequality)

$$
T_r(B^r)=T_{r,mod}(\Omega,u)\leq T_r(\Omega)\leq T_r(B_1) \left(\frac{|\Omega|}{\omega_N}\right)^{\frac{r(N+1) - N}{N}}.
$$
On the other hand, by scaling, calling $\rho(r)$ the radius of $B^r$, we have \[
T_r(B^r)=\rho(r)^{r(N+1)-N}T_r(B_1),
\]
where $B_1$ is the ball with unit radius centered in the origin. 
Therefore, we have that \[
\rho(r)\leq C(N) |\Omega|^{1/N},
\] 
and we conclude that $\rho(r)$ is uniformly bounded for $r\in (1,2]$, whence the balls $B^r$ are uniformly bounded as desired, namely all contained in a big ball $B_R$ for some $R>1$.

Now, we observe that, using H\"older inequality and the fact that $u\in W^{1,r}(\Omega)$ for all $r\in (1,2]$,
\[
\int_{B_R} |\nabla u_r^*|\, \de x = \int_{B^r}|\nabla u^*_r|\,\de x\leq |B^r|^{1/r'}\left(\int_{B^r}|\nabla u^*_r|^r\,\de x\right)^{1/r}\leq |B_R|^{1/r'}\left(\int_\Omega |\nabla u|^r\,\de x\right)^{1/r}		\leq C,
\]
for some positive constant $C=C(u,\Omega)$ independent of $r$. Thus $(u_r^{*})_r$ is uniformly bounded in $W^{1,1}_0(B_R)\subset BV_0(B_R)$, hence $u_r^*\rightharpoonup u^{\sharp}$ weakly in $ BV_0(B_R)$ as $r\to 1^+$ for some $u^\sharp\in BV_0(B_R)$. As a consequence there holds, using also the lower semicontinuity of the total variation,
		\begin{align}
			\|Du^{\sharp}\|_{TV(B_\sharp)}
			&\leq \liminf_{r\to 1}\int_{B^r}|\nabla u_r^*|\,\de x= \liminf_{r\to 1}\lim_{s\to 1}\int_{B^r}|\nabla u_r^*|^s\,\de x \notag\\ &=\liminf_{r\to 1}\lim_{s\to 1}\int_{\Omega}|\nabla u|^s\,\de x=\int_\Omega |\nabla u|\,\de x,\label{A}
			\\
			\int_{B_\sharp} |u^{\sharp}|^q\,\de x &= \liminf_{r\to1}\int_{B^r} |u_r^*|^q\,\de x\geq \int_\Omega |u|^q\,\de x,\label{B}
\end{align}
		for all $q\in[1,1^*)$, and calling $B_\sharp:=\{u^{\sharp}>0\}$. Notice that the support of $u^\sharp$, $B_\sharp$, is a ball.
		We are now in position to state our hypothesis.
		\begin{hypothesis}\label{mainhp}  Let $u\in \mathcal A(\Omega)$. Let $u^\sharp$ be the function constructed above and let $B_\sharp:=\{u^\sharp>0\}$. Let $B_{\Omega,u}$ be the ball such that $T_2(B_{\Omega,u})=T_{2,mod}(\Omega,u)$.
		
				We assume that
			\begin{equation}\label{eq:hp}
				B_\sharp\subseteq B_{\Omega,u}.
			\end{equation}
		\end{hypothesis}
		
		We shall apply Hypothesis \ref{mainhp}  to a particular family $(\varphi_n)_n$ of functions in $\mathcal A(\Omega)$, which are a minimizing sequence for the Cheeger constant of $\Omega$, that is, 
		\[
		\frac{\int_{\Omega}|\nabla \varphi_n|\, \de x}{\int_{\Omega}|\varphi_n|\,\de x}\longrightarrow \la_{1,1}(\Om),\qquad \mbox{as }n\rightarrow \infty\,.
		\] 
		The existence of such an approximating sequence $(\varphi_n)_n$ is not guaranteed in general, but it can be found for $\Omega$ regular enough and self--Cheeger, as we prove in the next Lemma.
			
		\begin{lemma}\label{apprAp}
			Let $\Omega$ be an open, bounded, convex set and self--Cheeger. For $n\in \N\setminus \{0\}$ we call $			\Omega_n:=\left\{x\in\Omega: {\rm dist}(x,\partial\Omega)>1/n  \right\},$ and we define  
			\[
			\varphi_n(x):=
			\begin{cases}
				1\quad&\text{ if }x\in\Omega_n\\
				n\,{\rm dist}(x,\partial\Omega)&\text{ if } x\in\Omega\setminus\Omega_n.
			\end{cases}
			\]
			For all $n\in\N$ such that $\Om_n\not=\emptyset$, we have $\varphi_n\in\mathcal A(\Omega)$  and moreover, as $n\to\infty$,
			\[
			\mathcal R(\varphi_n):=\frac{\int_{\Omega}|\nabla \varphi_n|\, \de x}{\int_{\Omega}|\varphi_n|\,\de x}\longrightarrow \la_{1,1}(\Omega).
			\]
		\end{lemma}	
		
		\begin{proof}
			The fact that $\varphi_n\in W^{1,2}_0(\Omega)\cap L^\infty(\Omega)$ is immediate because the functions are by construction Lipschitz continuous on a bounded set. Moreover, since
			\[
			|\nabla \varphi_n(x)|=
			\begin{cases}
				0\quad&\text{ if }x\in \Omega_n\\
				n&\text{ if }x\in \Omega\setminus \Om_n,
			\end{cases}
			\]
			one easily deduces that $\varphi_n\in \mathcal A_p(\Omega)$ for all $1<p\le 2$, because for all $t\in(0,1)$, the isoperimetric inequality entails that
			$$\frac{|\{x\in\Om\;:\;\varphi_n(x)>t\}|}{\int_{\{\varphi_n=t\}}|\nabla \varphi_n|^{p-1}\,\de \mathcal{H}^{N-1}}=\frac{|\{\varphi_n>t\}|}{n^{p-1}\,\mathcal H^{N-1}(\{\varphi_n=t\})}\leq \frac{C(N)}{n^{p-1}}|\{\varphi_n>t\}|^{1/N},
			$$
			so condition~\eqref{cattiva} is satisfied.
			We compute now
			\[
			\begin{aligned}
				\int_{\Omega}|\varphi_n|\, \de x&=	\int_{\Omega_n}1\,\de x+	\int_{\Omega\setminus\Omega_n}|\varphi_n| \de x =|\Omega_n|+\int_{\Omega\setminus\Omega_n}|\varphi_n|\,\de x \le |\Om_n|+|\Omega\setminus\Omega_n|\longrightarrow|\Om|\qquad \mbox{as }n\to \infty\,,
			\end{aligned}
			\]
			using also the fact that $\varphi_n\le1$. We observe moreover that, by coarea (and recalling that $\Omega_{1/t}=\{x\in \Omega : \dist(x,\partial\Omega)>t\}$),
			\[
			\begin{aligned}
				\int_{\Omega}|\nabla \varphi_n|\, \de x& =n\int_{\Omega\setminus \Omega_n}|\nabla \dist(x,\partial\Omega)|\, \de x=n\,|\Omega\setminus\Omega_n|
				=n\,\int_0^{1/n}\mathcal H^{N-1}(\partial\Omega_{1/t})\, \de t 
				\to\mathcal H^{N-1}(\partial\Omega)
			\end{aligned}
			\]
			as $n\to \infty$. Notice that the map $t\mapsto \mathcal H^{N-1}(\partial \Omega_{1/t})$ is continuous as $\Omega$ is convex. The conclusion of the Lemma easily follows from the assumption that $\Omega$ is self--Cheeger.
		\end{proof}

		\section{Proof of Theorem \ref{thm:ckjunderassumption}}\label{section6}
		
		This section is devoted to the proof of Theorem~\ref{thm:ckjunderassumption}, which shows the validity of Conjecture~\ref{conjecture}, under the assumption concerning the existence of a suitable rearrangement, explained in Section~\ref{sec-KJ}.
		
		\subsection{The Cheeger-Kohler-Jobin inequality for convex sets}
		
		We begin with a simple but useful result, which states that the minimization of $T_2(\cdot)^\theta\la_{1,1}(\cdot)$ among convex sets can be equivalently performed in the subclass of self--Cheeger sets, that is, if $\la_{1,1}(\Om)=\frac{P(\Om)}{|\Om|}$. 
		
		\begin{lemma}
			Let $\theta=\frac{1}{N+2}$. Then
			\begin{equation}\label{optselfcheeger}
				\begin{split}
					&\inf{\left\{T_2(\Om)^\theta\la_{1,1}(\Om)\;:\;\Om\subset\R^N,\;\mbox{open, bounded and convex}\right\}}\\
					&=\inf{\left\{T_2(\Om)^\theta\la_{1,1}(\Om)\;:\;\Om\subset\R^N,\;\mbox{open, bounded, convex and self--Cheeger}\right\}}\,.
				\end{split}
			\end{equation}
		\end{lemma}
		\begin{proof}
			The left hand side is clearly smaller than or equal to the right hand side. 
			In order to prove equality, we assume for the sake of contradiction that the strict inequality holds.  
			Then, given a convex set $\Om$, its Cheeger set $E_\Omega$ is convex too (see ~\cite{ac}) and by definition it satisfies $\la_{1,1}(E_\Om)=\la_{1,1}(\Om)$. 
			Moreover we can consider $E_\Om$ to be open without loss of generality, because the convexity entails $\la_{1,1}(E_\Om)=\la_{1,1}(E_\Om^{\circ})$; hence $E_\Om\subset\Om$.
			If $\Omega$ is not self--Cheeger, then $\cp_1(\Om\setminus E_\Om)>0$; in particular $\cp_2(\Om\setminus E_\Om)>0$, so that $T_2(\Om)>T_2(E_\Om)$. All in all, $E_\Omega$ is a better competitor for $ T_2(\cdot)^\theta \la_{1,1}(\cdot)$ in the left hand side in \eqref{optselfcheeger}. Clearly $E_\Om$ is self--Cheeger, thus the strict inequality does not hold and we have reached a contradiction.
		\end{proof}

		We are now in position to prove Theorem~\ref{thm:ckjunderassumption}, namely the Cheeger-Kohler-Jobin inequality, if Hypothesis~\ref{mainhp} holds.
\begin{proof}[Proof of Theorem \ref{thm:ckjunderassumption}] 
				We consider an open, bounded convex set $\Om\subset\R^N$. Without loss of generality we suppose it is also  self--Cheeger,  thanks to Lemma~\ref{optselfcheeger}. Thus, we can apply Lemma~\ref{apprAp}, in order to build a sequence of non-negative functions $(\varphi_n)_n \subset  \mathcal A(\Om)$, for $n\in \N$, such that
\[
\mathcal R(\varphi_n):=\frac{\int_{\Omega}|\nabla \varphi_n|\, \de x}{\int_{\Omega}|\varphi_n|\de x}\longrightarrow \la_{1,1}(\Om),\qquad \mbox{as }n\rightarrow \infty\,.
\]
				By applying the rearrangement described in Section~\ref{sec-KJ} to $\varphi_n\in \mathcal A(\Omega)$, exploiting the bounds \eqref{A} and \eqref{B} with $q=1$, we find a function $\varphi_n^\sharp$ such that
				$$
				\mathcal R(\varphi_n) \geq \mathcal R(\varphi_n^\sharp) = \frac{\int_{B_{\sharp,n}}|\nabla \varphi_n^\sharp|\, \de x}{\int_{B_{\sharp,n}}|\varphi_n^\sharp|\de x},
				$$
				where, with a slight abuse of notation, we have denoted $B_{\sharp,n}$ the ball $\{\varphi_n^\sharp>0\}$, while we will call $B_n$ the ball such that $T_2(B_n)=T_{2,mod}(\Om,\varphi_n)$ .
This yields, for $\theta=\frac{1}{N+2}$,
\[
					T_2(\Omega)^\theta\mathcal R(\varphi_n)\ge T_{2,mod}(\Omega, \varphi_n )^\theta\mathcal R(\varphi_n) =T_2(B_n)^\theta\mathcal R(\varphi_n)\geq T_2(B_n)^\theta\mathcal R(\varphi_n^\sharp) \geq T_2(B_n)^\theta \la_{1,1}(B_{\sharp,n})\,,
\]
				where the first inequality follows from Lemma~\ref{boiade}, while the last one is true since  $\varphi_n^\sharp$  is an admissible function for the infimum defining $\lambda_{1,1}(B_{\sharp,n})$. Using now Hypothesis~\ref{mainhp}, namely $B_{\sharp,n}\subset B_n$, and the monotonicity of $\la_{1,1}$, we deduce
				\begin{equation}\label{eq:chain2}
	T_2(\Omega)^\theta\mathcal R(\varphi_n)\geq T_2(B_n)^\theta \la_{1,1}(B_{\sharp,n})\ge  T_2(B_n)^\theta \la_{1,1}(B_n)\,.			
\end{equation}
				
				\noindent We observe that the quantity on the right-hand side of the chain of inequalities in \eqref{eq:chain2} is constant since the functional $T_2^\theta(\cdot)\la_{1,1}(\cdot)$ is scale invariant and $B_n$ are all balls.
				Hence, passing to the limit as $n\rightarrow\infty$ on the left-hand side, we obtain\[
				T_2(\Om)^\theta\la_{1,1}(\Om)=\limsup_{n\rightarrow \infty}T_2(\Omega)^\theta\mathcal R(\varphi_n)\geq T_2(B_n)^\theta\la_{1,1}(B_n)=T_2(B)^\theta\la_{1,1}(B)\,,
				\]
where $B$ is any ball.
\end{proof}

\section{Proof of the quantitative estimate for the Cheeger constant}\label{sec-qe}
		We offer here  the proof of Corollary~\ref{thm1}, assuming that Conjecture~\ref{conjecture} holds. We remark that it is a  modification of the combination of a Kohler-Jobin type inequality with a quantitative Saint-Venant inequality  proposed in~\cite{bdpv}.
		We recall that the quantitative Saint-Venant inequality proved in \cite[Section~5, Proof of Main Theorem]{bdpv} reads as
		\begin{equation}\label{qsv}
		T_2(B)|B|^{-\frac{N+2}{N}} -T_2(\Omega)|\Omega|^{-\frac{N+2}{N}} \ge \tau(N) \alpha(\Omega)^2\,,
		\end{equation}
		where $\tau=\tau(N)>0$ is a dimensional constant and $\alpha$ the Fraenkel asymmetry defined in~\eqref{fraenky}.

		\begin{proof}[Proof of Corollary  \ref{thm1}] 
			Since inequality~\eqref{qscheeger} is scale invariant, we may assume without loss of generality that $\Omega$ and $B$ have the same measure, say 1. Moreover, for brevity of notation, we will denote by $\theta$ the homogeneity exponent $\theta:=1/(N+2)$, and by $\tau$ the dimensional constant appearing in \eqref{qsv}. Thanks to the Cheeger--Kohler-Jobin estimate~\eqref{intro-KJq}, we have
			
			$$
			\frac{h_1(\Omega)}{h_1(B)} - 1 \geq \left (\frac{T_2(B)}{T_2(\Omega)} \right)^{\theta} - 1\,. 
			$$
			
			We now distinguish two cases: $T_2(B)/T_2(\Omega)> 2$ and  $T_2(B)/T_2(\Omega) \in [1,2]$.
			In the former, exploiting the easy bound $\tau \alpha^2(\Omega) \leq T_2(B) $, we obtain
			\begin{equation}\label{case1}
				\frac{h_1(\Omega)}{h_1(B)} - 1 \geq 2^{\theta} - 1   \geq  (2^{\theta}-1) \frac{\tau  \alpha^2(\Omega) }{T_2(B)}\,.
			\end{equation}
			In the latter, we use the concavity of the function $x\mapsto x^\theta$, being $0 < \theta=1/(N+2) < 1$.  For every $x\in [1,2]$, we have
			\begin{equation}\label{concave}
				x^{\theta} = (  (2-x) +  2 (x-1)  )^{\theta} \geq (2-x) + 2^{\theta}(x-1)\,,
			\end{equation}
			since $2-x$ and $x-1$ are both in $[0,1]$ and their sum is 1. By applying \eqref{concave} to $x={T_2(B)}/{T_2(\Omega)}$, we obtain
			\begin{equation}\label{case2}
				\begin{split}
					\frac{h_1(\Omega)}{h_1(B)} - 1 & \geq \left (2 - \frac{T_2(B)}{T_2(\Omega)}\right ) + 2^{\theta} \left( \frac{T_2(B)}{T_2(\Omega)} - 1 \right) - 1 = (2^{\theta}-1) \left ( \frac{T_2(B)}{T_2(\Omega)} - 1\right) 
					\\ & \geq (2^{\theta}-1) \frac{\tau \alpha^2(\Omega)}{T_2(\Omega)}\geq  (2^{\theta}-1) \frac{\tau \alpha^2(\Omega)}{T_2(B)}\,.
				\end{split}
			\end{equation}
			Finally, combining \eqref{case1} and \eqref{case2}, we conclude the proof of \eqref{qscheeger} with
			$$
			\sigma:= \frac{\tau  (2^{\theta}-1) h_1(B)}{T_2(B)}\,,
			$$
			where $B$ denotes an $N$-dimensional ball of unit measure.

			\noindent Eventually, we notice that the exponent $2$ is sharp.  Indeed it is enough to consider the family of ellipsoids 
			\[
			\Omega_\eps=\gamma_\eps\Big\{(x_1,\dots,x_N): \sum_{i=1}^{N-1}x_i^2+(1+\eps)x_N^2\le 1\Big\}\,,
			\]
			where $\gamma_\eps>0$ is such that $|\Omega_\eps|=1$.
			A simple computation shows that $P(\Omega_\eps)-P(B)\simeq \eps^2$ while $\alpha(\Omega_\eps)\simeq \eps$ as $\eps\to 0$. On the other hand, since $B$ is self-Cheeger,
			\[
			h_1(\Omega_\eps)-h_1(B)\leq P(\Omega_\eps)-P(B)\simeq \eps^2\,.
			\] 
			This concludes the proof. 
		\end{proof}

	\end{document}